\documentclass[reqno,a4paper]{amsart}
\usepackage{amssymb}%
\usepackage[hmarginratio=1:1]{geometry}%
\usepackage{ifpdf}
\usepackage{listings}
\usepackage{color,xcolor}
\usepackage{colortbl}
\usepackage{multirow}
\usepackage{booktabs}
\usepackage{float}
\usepackage[colorlinks,linkcolor=red,anchorcolor=blue,citecolor=red]{hyperref}
\usepackage[numbers,sort&compress]{natbib}
\allowdisplaybreaks[4]
\theoremstyle{plain}

\newtheorem{cor}{Corollary}[section]

\newtheorem{mr}{Monotonicity rule}[section]
\theoremstyle{remark}
\newtheorem{rem}{Remark}[section]
\theoremstyle{plain}

\theoremstyle{definition}

\numberwithin{equation}{section}

\begin{document}

\title[]{Monotonicity rules for the ratio of power series}
\author{Zhong-Xuan Mao, Jing-Feng Tian*}

\address{Zhong-Xuan Mao \\
Hebei Key Laboratory of Physics and Energy Technology\\
Department of Mathematics and Physics\\
North China Electric Power University \\
Yonghua Street 619, 071003, Baoding, P. R. China}
\email{maozhongxuan000\symbol{64}gmail.com}

\address{Jing-Feng Tian\\
Hebei Key Laboratory of Physics and Energy Technology\\
Department of Mathematics and Physics\\
North China Electric Power University\\
Yonghua Street 619, 071003, Baoding, P. R. China}
\email{tianjf\symbol{64}ncepu.edu.cn}

\begin{abstract}
In this paper, we present some monotonicity rules for the ratio of two power series $x\mapsto \sum_{k=0}^\infty a_k x^k / \sum_{k=0}^\infty b_k x^k$ under the assumption that the monotonicity of the sequence ${a_k/b_k}$ changes twice. Additionally, we introduce a local monotonicity rule in this paper.
\end{abstract}


\thanks{*Corresponding author: Jing-Feng Tian(tianjf\symbol{64}ncepu.edu.cn)}
\maketitle

\section{Introduction}
Monotonicity rules are essential in the field of analytics, and widely used in fields such as approximation theory, differential geometry, information theory, probability, and statistics. Meanwhile, they are crucial tools for exploring the properties of special functions.

The origins of the monotonicity rules for the ratio of functions can be attributed to a lemma
\cite[Lemma 1]{Biernacki-AUMCS-1955} that was initially employed by Biernacki and Krzy\.z in their studies on differential geometry.
\begin{mr} \label{mr-SER-1}
Let real power series $\mathcal{A}(x)=\sum_{k=0}^\infty a_k x^k$ and $\mathcal{B}(x)=\sum_{k=0}^\infty b_k x^k$ converge on $(-r,r)$, and $b_k$ be non-negative and not vanish identically. If the sequence $\{a_k/b_k\}$ is increasing, then the function
\begin{equation} \label{mr-SER-1-f-1}
x \mapsto \frac{\mathcal{A}(x)}{\mathcal{B}(x)} = \frac{\sum_{k=0}^\infty a_k x^k }{\sum_{k=0}^\infty b_k x^k}
\end{equation}
is increasing on $(0,r)$.
\end{mr}

Next we introduce the following monotonicity rule (see \cite{Yang-JMAA-2015} and \cite{Yang-RJ-2019}), which considers the case that the monotonicity of $\{a_k/b_k\}_{k\geq0}$ changes once on the basis of Monotonicity rule \ref{mr-SER-1}.
\begin{mr} \label{mr-SER-5}
Let real power series $\mathcal{A}(x)=\sum_{k=0}^\infty a_{k} x^{k}$ and $\mathcal{B}(x)=\sum_{k=0}^\infty b_{k} x^{k}$ be converge on $(-r,r)$ with $b_k>0$. If there exists a given integer $m\geq1$ such that $\{a_k/b_k\}$ is increasing (decreasing) for all $0\leq k\leq m$ and decreasing (increasing) for all $k\geq m$ with neither  $\{a_k/b_k\}_{0\leq k \leq m}$ nor $\{a_k/b_k\}_{k \geq m}$ is constant, then
\begin{itemize}
\item[(1)] the function $x \mapsto \mathcal{A}(x)/\mathcal{B}(x)$ is increasing (decreasing) on $(0,r)$ if and only if $H_{\mathcal{A},\mathcal{B}}(r^-) \geq (\leq) 0$,
\item[(2)] there exists $x_0\in(0,r)$ such that the function $x \mapsto \mathcal{A}(x)/\mathcal{B}(x)$ is increasing (decreasing) on $(0,x_0]$ and decreasing (increasing) on $[x_0,r)$ if $H_{\mathcal{A},\mathcal{B}}(r^-) < (>) 0$,
\end{itemize}
where the function $H_{f,g}$ is defined as $H_{f,g}:=\frac{f^\prime}{g^\prime}g-f$, which is called Yang's $H$-function (named in \cite{Tian-AIMS-2020}).
\end{mr}

\begin{rem}
Monotonicity rule \ref{mr-SER-5} is first introduced by Yang, Chu, and Wang in \cite[Theorem 2.1]{Yang-JMAA-2015}, and it is slightly modified by Yang and Tian in \cite[Lemma 2]{Yang-RJ-2019}. Here we give the modified version.
\end{rem}

In the perspective of the number of times the monotonicity of $\{a_k/b_k\}$ undergoes changes, a research gap becomes evident within the existing research: concerning scenarios where such monotonicity changes occur 2 times.
The primary objective of this paper is to bridge the existing gap and extend the applicability of the monotonicity rules for the ratio of power series.

In this paper, we will establish the monotonicity rules for
\begin{equation} \label{AB}
x \mapsto \frac{\mathcal{A}(x)}{\mathcal{B}(x)} = \frac{\sum_{k=0}^\infty a_k x^k }{\sum_{k=0}^\infty b_k x^k},
\end{equation}
under the case that the monotonicity of $\{a_k/b_k\}$ change twice. We also provide a local monotonicity rules in this section.

\section{Monotonicity rules for the ratio of power series}

We begin with a complete introduction to the Yang's $H$ function \cite{Yang-JMAA-2015}. Let $-\infty \leq a < b\leq \infty$, functions $F$ and $G$ be differentiable on $(a,b)$, and $G\neq0$ on $(a,b)$. Yang's $H$ function is defined by
\begin{equation*}
H_{F,G} =\frac{F^\prime}{G^\prime} G -F.
\end{equation*}
It is easy to check that the following two formulas hold:
\begin{equation*}
\Big( \frac{F}{G} \Big)^\prime = \frac{G^\prime}{G^2} H_{F,G},
\end{equation*}
and
\begin{equation*}
H_{F,G}^\prime= \Big( \frac{F^\prime}{G^\prime} \Big)^\prime G,
\end{equation*}
where the second identity requires both $F$ and $G$ are twice differentiable.
In what follows, we also denote that $[a,b]_\mathbb{N}:=[a,b] \cap \mathbb{N}$, where $\mathbb{N}:=\{0,1,2,\cdots\}$ is the set of natural numbers.
Here and after, $\mathcal{A}^\prime$ and $\mathcal{B}^\prime$ are the derivatives of functions $\mathcal{A}$ and $\mathcal{B}$, respectively. We also use the forms $\mathcal{A}^{(k)}$ and $\mathcal{B}^{(k)}$ to respectively represent the derivatives of $\mathcal{A}$ and $\mathcal{B}$, where $\mathcal{A}^{(0)}$ and $\mathcal{B}^{(0)}$ means $\mathcal{A}$ and $\mathcal{B}$.

\subsection{The monotonicity of $\{a_k/b_k\}$ changes twice}
First, we establish the monotonicity rule for function (\ref{AB}), where the monotonicity of $\{a_k/b_k\}$ changes twice.

\begin{mr} \label{DW21}
Let real power series $\mathcal{A}(x)=\sum_{k=0}^\infty a_k x^k$ and $\mathcal{B}(x)=\sum_{k=0}^\infty b_k x^k$ converge on $(0,r)$ with $b_k>0$. If there exists different integers $m_2>m_1\geq1$ such that the sequence $\{a_k/b_k\}$ is increasing (decreasing) for all $k \in [0,m_1]_\mathbb{N} \cup [m_2,\infty)_\mathbb{N}$ and decreasing (increasing) for all $k \in [m_1,m_2]_\mathbb{N}$, as well as $\{a_k/b_k\}_{0\leq k \leq m_1}$, $\{a_k/b_k\}_{m_1 \leq k \leq m_2}$, and $\{a_k/b_k\}_{k \geq m_2}$ are non-constant, then
\begin{itemize}
\item[(C1)] the function $x\mapsto \mathcal{A}(x)/\mathcal{B}(x)$ is increasing (decreasing) on $(0,r)$ if one of the following conditions holds:
    \begin{itemize}
    \item[(i)] $H_{\mathcal{A}, \mathcal{B}}(r^-) \geq (\leq) 0$ and  $H_{\mathcal{A}^\prime, \mathcal{B}^\prime}(r^-) \leq (\geq) 0$;
    \item[(ii)] $H_{\mathcal{A}, \mathcal{B}}(r^-) > (<) 0$, $H_{\mathcal{A}^\prime, \mathcal{B}^\prime}(r^-) > (<) 0$, and $H_{\mathcal{A}, \mathcal{B}}(x) \geq (\leq) 0$ for all $x\in(0,r)$.
    \end{itemize}
\item[(C2)] there exists $x_1 \in (0,r)$ such that the function $x\mapsto \mathcal{A}(x)/\mathcal{B}(x)$ is increasing (decreasing) on $(0,x_1]$ and decreasing (increasing) on $[x_1,r)$ if one of the following conditions holds:
    \begin{itemize}
    \item[(iii)] $H_{\mathcal{A}, \mathcal{B}}(r^-) < (>) 0$ and $H_{\mathcal{A}^\prime, \mathcal{B}^\prime}(r^-) \leq (\geq) 0$;
    \item[(iv)] $H_{\mathcal{A}, \mathcal{B}}(r^-) \leq (\geq) 0$ and $H_{\mathcal{A}^\prime, \mathcal{B}^\prime}(r^-) > (<) 0$.
    \end{itemize}
\item[(C3)] there exists $x_2,x_3 \in (0,r)$ such that the function $x\mapsto \mathcal{A}(x)/\mathcal{B}(x)$ is increasing (decreasing) on $(0,x_2] \cup [x_3,r)$ and decreasing (increasing) on $[x_2,x_3]$ if:
    \begin{itemize}
    \item[(v)] $H_{\mathcal{A}, \mathcal{B}}(r^-) \geq (\leq) 0$, $H_{\mathcal{A}^\prime, \mathcal{B}^\prime}(r^-) > (<) 0$, and there exists $x_0\in(0,r)$ such that $H_{\mathcal{A}, \mathcal{B}}(x_0) < (>) 0$.
    \end{itemize}
\end{itemize}
\end{mr}

\begin{proof}
Without loss of generality, we consider the case that $\{a_k/b_k\}$ is increasing for all $k \in [0,m_1]_\mathbb{N} \cup [m_2,\infty)_\mathbb{N}$ and decreasing for all $k \in [m_1,m_2]_\mathbb{N}$.
We complete the proof by using mathematical induction for integer $m_1$.

(I) When $m_1=1$, that is, the sequence $\{a_k/b_k\}$ is increasing for all $k \in \{0,1\} \cup [m_2,\infty)_\mathbb{N}$ and decreasing for all $k \in [1,m_2]_\mathbb{N}$ as well as the sequence $\{a_{k+1}/b_{k+1}\}$ is decreasing for all $0 \leq k \leq m_2-1$ and increasing for all $k \geq m_2-1$. Noting that
\begin{equation*}
\frac{\mathcal{A}^\prime(x)}{\mathcal{B}^\prime(x)} = \frac{\sum_{k=0}^\infty (k+1) a_{k+1} x^{k} }{\sum_{k=0}^\infty (k+1) b_{k+1} x^{k}},
\end{equation*}
by using Monotonicity rule \ref{mr-SER-5}, we obtain that
\begin{itemize}
\item[(i)] the function $x\mapsto \mathcal{A}^\prime(x)/\mathcal{B}^\prime(x)$ is decreasing on $x\in(0,r)$ if $H_{\mathcal{A}^\prime, \mathcal{B}^\prime}(r^-) \leq 0$,
\item[(ii)] there exists $x_0\in(0,r)$ such that it is decreasing on $(0,x_0]$ and increasing on $[x_0,r)$ if $H_{\mathcal{A}^\prime, \mathcal{B}^\prime}(r^-) > 0$.
\end{itemize}
From identity $H_{\mathcal{A},\mathcal{B}}^\prime= (\mathcal{A}^\prime/\mathcal{B}^\prime)^\prime \mathcal{B}$ and $\mathcal{B}>0$, we obtain that the monotonicity of $x\mapsto H_{\mathcal{A},\mathcal{B}}(x)$ is the same as $x\mapsto \mathcal{A}^\prime(x)/\mathcal{B}^\prime(x)$.
Noting that
\begin{equation*}
H_{\mathcal{A}, \mathcal{B}}(0^+) = \frac{\mathcal{A}^\prime(0^+)}{\mathcal{B}^\prime(0^+)}\mathcal{B}(0^+)- \mathcal{A}(0^+) = b_0 \Big( \frac{a_1}{b_1} - \frac{a_0}{b_0} \Big) > 0,
\end{equation*}
we receive the following five conclusions:

\begin{itemize}
\item[(i)] If $H_{\mathcal{A}, \mathcal{B}}(r^-) \geq 0$ and $H_{\mathcal{A}^\prime, \mathcal{B}^\prime}(r^-) \leq 0$, then $H_{\mathcal{A}, \mathcal{B}} \geq 0$. From identity $(\mathcal{A}/\mathcal{B})^\prime = \mathcal{B}^\prime/ \mathcal{B}^2 H_{\mathcal{A}, \mathcal{B}}$, we obtain that $x\mapsto \mathcal{A}(x)/\mathcal{B}(x)$ is increasing on $(0,r)$.
\item[(ii)] If $H_{\mathcal{A}, \mathcal{B}}(r^-) < 0$ and $H_{\mathcal{A}^\prime, \mathcal{B}^\prime}(r^-) \leq 0$, then there exists $x_1 \in (0,r)$ such that $H_{\mathcal{A}, \mathcal{B}}(x) \geq 0$ for all $x\in(0,x_1]$ and $H_{\mathcal{A}, \mathcal{B}}(x) \leq 0$ for all $x\in[x_1,r)$. Moreover, the function $x\mapsto \mathcal{A}(x)/\mathcal{B}(x)$ is increasing on $(0,x_1]$ and it is decreasing on $[x_1,r)$.
\item[(iii)] If $H_{\mathcal{A}, \mathcal{B}}(r^-) > 0$, $H_{\mathcal{A}^\prime, \mathcal{B}^\prime}(r^-) > 0$ and $H_{\mathcal{A}, \mathcal{B}}(x_0) \geq 0$, then $H_{\mathcal{A}, \mathcal{B}} \geq 0$. Moreover, the function $x\mapsto \mathcal{A}(x)/\mathcal{B}(x)$ is increasing on $(0,r)$.
\item[(iv)] If $H_{\mathcal{A}, \mathcal{B}}(r^-) > 0$, $H_{\mathcal{A}^\prime, \mathcal{B}^\prime}(r^-) > 0$ and $H_{\mathcal{A}, \mathcal{B}}(x_0) < 0$, then there exist $x_2,x_3\in(0,r)$ such that $H_{\mathcal{A}, \mathcal{B}}(x) \geq 0$ for all $x\in(0,x_2] \cup [x_3,r)$ and $H_{\mathcal{A}, \mathcal{B}}(x) \leq 0$ for all $x\in[x_2,x_3]$. Moreover, the function $x\mapsto \mathcal{A}(x)/\mathcal{B}(x)$ is increasing on $(0,x_2] \cup [x_3,r)$ and it is decreasing on $[x_2,x_3]$.
\item[(v)] If $H_{\mathcal{A}, \mathcal{B}}(r^-) \leq 0$ and $H_{\mathcal{A}^\prime, \mathcal{B}^\prime}(r^-) > 0$, then there exists $x_1 \in (0,r)$ such that $H_{\mathcal{A}, \mathcal{B}}(x) \geq 0$ for all $x\in(0,x_1]$ and $H_{\mathcal{A}, \mathcal{B}}(x) \leq 0$ for all $x\in[x_1,r)$. Moreover, the function $x\mapsto \mathcal{A}(x)/\mathcal{B}(x)$ is increasing on $(0,x_1]$ and it is decreasing on $[x_1,r)$.
\end{itemize}

(II) Suppose desired conclusions hold for $m_1=n$, namely, $\{a_k/b_k\}$ is increasing for all $k \in [0,n]_\mathbb{N} \cup [m_2,\infty)_\mathbb{N}$ and decreasing for all $k \in [n,m_2]_\mathbb{N}$, as well as (C1), (C2), and (C3) hold.

(III) When $m_1=n+1$, namely, $\{a_k/b_k\}$ is increasing for all $k \in [0,n+1]_\mathbb{N} \cup [m_2,\infty)_\mathbb{N}$ and decreasing for all $k \in [n+1,m_2]_\mathbb{N}$. In this case, the sequence $\{a_{k+1}/b_{k+1}\}$ is increasing for all $k \in [0,n]_\mathbb{N} \cup [m_2-1,\infty)_\mathbb{N}$ and decreasing for all $k \in [n,m_2-1]_\mathbb{N}$. Thus, by using induction hypothesis, we know that (C1), (C2), and (C3) hold for function $x\mapsto \mathcal{A}^\prime(x)/\mathcal{B}^\prime(x)$ and obtain the following conclusions:
\begin{itemize}
\item[(i)] If $H_{\mathcal{A}^\prime,\mathcal{B}^\prime}(r^-) < 0$, then there exists $x_1\in(0,r)$ such that the function $x\mapsto \mathcal{A}^\prime(x)/\mathcal{B}^\prime(x)$ is increasing on $(0,x_1]$ and decreasing on $[x_1,r)$, same monotonicity as $H_{\mathcal{A},\mathcal{B}}$. Thus, $H_{\mathcal{A},\mathcal{B}}\geq0$ if $H_{\mathcal{A},\mathcal{B}}(r^-)\geq0$, as well as $H_{\mathcal{A},\mathcal{B}}(x)\geq0$ for all $x\in (0,x_1]$ and $H_{\mathcal{A},\mathcal{B}}(x)\leq0$ for all $x\in[x_1,r)$ if $H_{\mathcal{A},\mathcal{B}}(r^-)<0$, namely, the function $x\mapsto \mathcal{A}(x)/\mathcal{B}(x)$ is increasing on $(0,r)$ if $H_{\mathcal{A},\mathcal{B}}(r^-)\geq0$ as well as it is increasing on $(0,x_1]$ and decreasing on $[x_1,r)$ if $H_{\mathcal{A},\mathcal{B}}(r^-)<0$.
\item[(ii)] If $H_{\mathcal{A}^\prime,\mathcal{B}^\prime}(r^-) > 0$, then the function $x\mapsto \mathcal{A}^\prime(x)/\mathcal{B}^\prime(x)$ is increasing on $(0,r)$ or there exists $x_2,x_3\in(0,r)$ such that it is increasing $(0,x_2] \cup [x_3,r)$ and decreasing on $[x_2,x_3]$.

    For the first case, the function $x\mapsto \mathcal{A}(x)/\mathcal{B}(x)$ is increasing on $(0,r)$. For the second case, if there exists $x_0\in(0,r)$ such that $H_{\mathcal{A},\mathcal{B}}(x_0) < 0$, then there exist $x_2,x_3\in(0,r)$ such that the function $x\mapsto \mathcal{A}(x)/\mathcal{B}(x)$ is increasing on $(0,x_2] \cup [x_3,r)$ and decreasing on $[x_2,x_3]$ if $H_{\mathcal{A},\mathcal{B}}(r^-)>0$, and there exist $x_1 \in (0,r)$ such that the function $x\mapsto \mathcal{A}(x)/\mathcal{B}(x)$ is increasing on $(0,x_1]$ and decreasing on $[x_1,r)$ if $H_{\mathcal{A},\mathcal{B}}(r^-)\leq0$. If $H_{\mathcal{A},\mathcal{B}}(x) \geq 0$ for all $x\in(0,r)$, then the function $x\mapsto \mathcal{A}(x)/\mathcal{B}(x)$ is increasing on $(0,r)$.
\item[(iii)] If $H_{\mathcal{A}^\prime,\mathcal{B}^\prime}(r^-) =0$, then the function $x\mapsto \mathcal{A}^\prime(x)/\mathcal{B}^\prime(x)$ is increasing on $(0,r)$ or there exists $x_1\in(0,r)$ such that the function $x\mapsto \mathcal{A}^\prime(x)/\mathcal{B}^\prime(x)$ is increasing on $(0,x_1]$ and decreasing $[x_1,r)$. Thus, the function $x\mapsto \mathcal{A}(x)/\mathcal{B}(x)$ is increasing on $(0,r)$ if $H_{\mathcal{A},\mathcal{B}}(r^-) \geq0$, and there exists $x_1\in(0,r)$ such that the function $x\mapsto \mathcal{A}(x)/\mathcal{B}(x)$ is increasing on $(0,x_1]$ and decreasing $[x_1,r)$ if $H_{\mathcal{A},\mathcal{B}}(r^-) < 0$.
\end{itemize}
Thus, we complete the proof.
\end{proof}

\begin{rem}
If one of $\{a_k/b_k\}_{0 \leq k \leq m_1}$, $\{a_k/b_k\}_{m_1 \leq k \leq m_2}$, and $\{a_k/b_k\}_{m_2 \leq k}$ is a constant, then we regard the monotonicity of $\{a_k/b_k\}$ change once rather than twice.
\end{rem}

\begin{rem}
Noting that the condition $H_{\mathcal{A}, \mathcal{B}}(r^-) = 0$ and $H_{\mathcal{A}^\prime, \mathcal{B}^\prime}(r^-) > 0$ deduces there exist $x_0 \in (0,r)$ such that $H_{\mathcal{A}, \mathcal{B}}(x_0) < 0$, we know that, for any given functions $\mathcal{A}$ and $\mathcal{B}$, one of (i), (ii), (iii), (iv), (v) in (C1), (C2), (C3) holds.
\end{rem}

\begin{rem}
Let
\begin{itemize}
\item[(i)] arrows ``$\nearrow$'' and ``$\searrow$'' respectively denote ``increasing'' and ``decreasing'';
\item[(ii)] arrows ``$\nearrow\searrow$'' and ``$\searrow\nearrow$'' respectively denote ``there exist $x_1\in(0,r)$ such that the function increasing on $(0,x_1]$ and decreasing on $[x_1,r)$'' and ``there exist $x_1\in(0,r)$ such that the function decreasing on $(0,x_1]$ and increasing on $[x_1,r)$'';
\item[(iii)] arrows ``$\nearrow\searrow\nearrow$'' and ``$\searrow\nearrow\searrow$'' respectively denote ``there exist $x_2,x_3\in(0,r)$ such that the function increasing on $(0,x_2]\cup [x_3,r)$ and decreasing on $[x_2,x_3]$'' and ``there exist $x_2,x_3\in(0,r)$ such that the function decreasing on $(0,x_2]\cup [x_3,r)$ and increasing on $[x_2,x_3]$''.
\end{itemize}
Then the conclusions in Monotonicity rule \ref{DW21} could be simply represented in Table \ref{table1}.
\begin{table}[h]

\begin{tabular}{cccccc}
  \arrayrulecolor{blue} \toprule[0.8pt]
  Cases & $\{a_k/b_k\}$ & $H_{\mathcal{A},\mathcal{B}}(r^-)$ & $H_{\mathcal{A}^\prime,\mathcal{B}^\prime}(r^-)$ & $H_{\mathcal{A},\mathcal{B}}(x)$ & $\mathcal{A}/\mathcal{B}$ \\
  \arrayrulecolor{blue} \midrule[0.8pt]
  $1$ & $\nearrow\searrow\nearrow$ & $\geq 0$ & $\leq 0$ & & $\nearrow$ \\
  $2$ & $\nearrow\searrow\nearrow$ & $> 0$    & $>0$     & $\geq 0$ for all $x\in(0,r)$ & $\nearrow$  \\
  $3$ & $\nearrow\searrow\nearrow$ & $< 0$    & $\leq 0$ & & $\nearrow\searrow$ \\
  $4$ & $\nearrow\searrow\nearrow$ & $\leq 0$ & $> 0$    & & $\nearrow\searrow$ \\
  $5$ & $\nearrow\searrow\nearrow$ & $\geq 0$ & $> 0$    & exists $x_0\in(0,r)$ such that $< 0$ & $\nearrow\searrow\nearrow$ \\
  \hline
  $6$ & $\searrow\nearrow\searrow$ & $\leq 0$ & $\geq 0$ & & $\searrow$ \\
  $7$ & $\searrow\nearrow\searrow$ & $< 0$    & $<0$     & $\leq 0$ for all $x\in(0,r)$ & $\searrow$  \\
  $8$ & $\searrow\nearrow\searrow$ & $> 0$    & $\geq 0$ & & $\searrow\nearrow$ \\
  $9$ & $\searrow\nearrow\searrow$ & $\geq 0$ & $< 0$    & & $\searrow\nearrow$ \\
  $10$ & $\searrow\nearrow\searrow$ & $\leq 0$ & $< 0$    & exists $x_0\in(0,r)$ such that $> 0$ & $\searrow\nearrow\searrow$ \\
  \arrayrulecolor{blue} \bottomrule[0.8pt]
\end{tabular}
\caption{The Monotonicity of $\mathcal{A}/\mathcal{B}$ in Monotonicity rule \ref{DW21}} \label{table1}
\end{table}
\end{rem}

Let $r\to \infty$. Then we have the following monotonicity rule.
\begin{mr}
Let real power series $\mathcal{A}(x)=\sum_{k=0}^\infty a_k x^k$ and $\mathcal{B}(x)=\sum_{k=0}^\infty b_k x^k$ converge on $(0,\infty)$ with $b_k>0$. If there exists different integers $m_2>m_1\geq1$ such that the sequence $\{a_k/b_k\}$ is increasing (decreasing) for all $k \in [0,m_1]_\mathbb{N} \cup [m_2,\infty)_\mathbb{N}$ and decreasing (increasing) for all $k \in [m_1,m_2]_\mathbb{N}$, as well as $\{a_k/b_k\}_{0\leq k \leq m_1}$, $\{a_k/b_k\}_{m_1 \leq k \leq m_2}$, and $\{a_k/b_k\}_{k \geq m_2}$ are non-constant., then
\begin{itemize}
\item[(C1)] the function $x\mapsto \mathcal{A}(x)/\mathcal{B}(x)$ is increasing (decreasing) on $(0,\infty)$ if one of the following conditions holds:
    \begin{itemize}
    \item[(i)] $H_{\mathcal{A}, \mathcal{B}}(\infty) \geq (\leq) 0$ and  $H_{\mathcal{A}^\prime, \mathcal{B}^\prime}(\infty) \leq (\geq) 0$;
    \item[(ii)] $H_{\mathcal{A}, \mathcal{B}}(\infty) > (<) 0$, $H_{\mathcal{A}^\prime, \mathcal{B}^\prime}(\infty) > (<) 0$, and $H_{\mathcal{A}, \mathcal{B}}(x) \geq (\leq) 0$ for all $x\in(0,\infty)$.
    \end{itemize}
\item[(C2)] there exists $x_1 \in (0,\infty)$ such that the function $x\mapsto \mathcal{A}(x)/\mathcal{B}(x)$ is increasing (decreasing) on $(0,x_1]$ and decreasing (increasing) on $[x_1,\infty)$ if one of the following conditions holds:
    \begin{itemize}
    \item[(iii)] $H_{\mathcal{A}, \mathcal{B}}(\infty) < (>) 0$ and $H_{\mathcal{A}^\prime, \mathcal{B}^\prime}(\infty) \leq (\geq) 0$;
    \item[(iv)] $H_{\mathcal{A}, \mathcal{B}}(\infty) \leq (\geq) 0$ and $H_{\mathcal{A}^\prime, \mathcal{B}^\prime}(\infty) > (<) 0$.
    \end{itemize}
\item[(C3)] there exists $x_2,x_3 \in (0,\infty)$ such that the function $x\mapsto \mathcal{A}(x)/\mathcal{B}(x)$ is increasing (decreasing) on $(0,x_2] \cup [x_3,r)$ and decreasing (increasing) on $[x_2,x_3]$ if:
    \begin{itemize}
    \item[(v)] $H_{\mathcal{A}, \mathcal{B}}(\infty) \geq (\leq) 0$, $H_{\mathcal{A}^\prime, \mathcal{B}^\prime}(\infty) > (<) 0$, and there exists $x_0\in(0,\infty)$ such that $H_{\mathcal{A}, \mathcal{B}}(x_0) < (>) 0$.
    \end{itemize}
\end{itemize}
\end{mr}

Likewise, we present the monotonicity rule for the ratio of two polynomials.
\begin{mr}
Let $\mathcal{A}_N(x)=\sum_{k=0}^N a_k x^k$ and $\mathcal{B}_N(x)=\sum_{k=0}^N b_k x^k$ defined on $(0,r)$ with $b_k>0$. If there exists different integers $N> m_2>m_1\geq1$ such that the sequence $\{a_k/b_k\}$ is increasing (decreasing) for all $k \in [0,m_1]_\mathbb{N} \cup [m_2,N]_\mathbb{N}$ and decreasing (increasing) for all $k \in [m_1,m_2]_\mathbb{N}$, as well as $\{a_k/b_k\}_{0\leq k \leq m_1}$, $\{a_k/b_k\}_{m_1 \leq k \leq m_2}$, and $\{a_k/b_k\}_{ m_2 \leq k \leq N}$ are non-constant., then
\begin{itemize}
\item[(C1)] the function $x\mapsto \mathcal{A}_N(x)/\mathcal{B}_N(x)$ is increasing (decreasing) on $(0,r)$ if one of the following conditions holds:
    \begin{itemize}
    \item[(i)] $H_{\mathcal{A}_N, \mathcal{B}_N}(r^-) \geq (\leq) 0$ and  $H_{\mathcal{A}_N^\prime, \mathcal{B}_N^\prime}(r^-) \leq (\geq) 0$;
    \item[(ii)] $H_{\mathcal{A}_N, \mathcal{B}_N}(r^-) > (<) 0$, $H_{\mathcal{A}_N^\prime, \mathcal{B}_N^\prime}(r^-) > (<) 0$, and $H_{\mathcal{A}_N, \mathcal{B}_N}(x) \geq (\leq) 0$ for all $x\in(0,r)$.
    \end{itemize}
\item[(C2)] there exists $x_1 \in (0,r)$ such that the function $x\mapsto \mathcal{A}_N(x)/\mathcal{B}_N(x)$ is increasing (decreasing) on $(0,x_1]$ and decreasing (increasing) on $[x_1,r)$ if one of the following conditions holds:
    \begin{itemize}
    \item[(iii)] $H_{\mathcal{A}_N, \mathcal{B}_N}(r^-) < (>) 0$ and $H_{\mathcal{A}_N^\prime, \mathcal{B}_N^\prime}(r^-) \leq (\geq) 0$;
    \item[(iv)] $H_{\mathcal{A}_N, \mathcal{B}_N}(r^-) \leq (\geq) 0$ and $H_{\mathcal{A}_N^\prime, \mathcal{B}_N^\prime}(r^-) > (<) 0$.
    \end{itemize}
\item[(C3)] there exists $x_2,x_3 \in (0,r)$ such that the function $x\mapsto \mathcal{A}_N(x)/\mathcal{B}_N(x)$ is increasing (decreasing) on $(0,x_2] \cup [x_3,r)$ and decreasing (increasing) on $[x_2,x_3]$ if:
    \begin{itemize}
    \item[(v)] $H_{\mathcal{A}_N, \mathcal{B}_N}(r^-) \geq (\leq) 0$, $H_{\mathcal{A}_N^\prime, \mathcal{B}_N^\prime}(r^-) > (<) 0$, and there exists $x_0\in(0,r)$ such that $H_{\mathcal{A}_N, \mathcal{B}_N}(x_0) < (>) 0$.
    \end{itemize}
\end{itemize}
\end{mr}

\begin{mr}
Let $\mathcal{A}_N(x)=\sum_{k=0}^N a_k x^k$ and $\mathcal{B}_N(x)=\sum_{k=0}^N b_k x^k$ defined on $(0,\infty)$ with $b_k>0$. If there exists different integers $N> m_2>m_1\geq1$ such that the sequence $\{a_k/b_k\}$ is increasing (decreasing) for all $k \in [0,m_1]_\mathbb{N} \cup [m_2,N]_\mathbb{N}$ and decreasing (increasing) for all $k \in [m_1,m_2]_\mathbb{N}$, as well as $\{a_k/b_k\}_{0\leq k \leq m_1}$, $\{a_k/b_k\}_{m_1 \leq k \leq m_2}$, and $\{a_k/b_k\}_{ m_2 \leq k \leq N}$ are non-constant., then
\begin{itemize}
\item[(C1)] the function $x\mapsto \mathcal{A}_N(x)/\mathcal{B}_N(x)$ is increasing (decreasing) on $(0,\infty)$ if one of the following conditions holds:
    \begin{itemize}
    \item[(i)] $H_{\mathcal{A}_N, \mathcal{B}_N}(\infty) \geq (\leq) 0$ and  $H_{\mathcal{A}_N^\prime, \mathcal{B}_N^\prime}(\infty) \leq (\geq) 0$;
    \item[(ii)] $H_{\mathcal{A}_N, \mathcal{B}_N}(\infty) > (<) 0$, $H_{\mathcal{A}_N^\prime, \mathcal{B}_N^\prime}(\infty) > (<) 0$, and $H_{\mathcal{A}_N, \mathcal{B}_N}(x) \geq (\leq) 0$ for all $x\in(0,\infty)$.
    \end{itemize}
\item[(C2)] there exists $x_1 \in (0,\infty)$ such that the function $x\mapsto \mathcal{A}_N(x)/\mathcal{B}_N(x)$ is increasing (decreasing) on $(0,x_1]$ and decreasing (increasing) on $[x_1,\infty)$ if one of the following conditions holds:
    \begin{itemize}
    \item[(iii)] $H_{\mathcal{A}_N, \mathcal{B}_N}(\infty) < (>) 0$ and $H_{\mathcal{A}_N^\prime, \mathcal{B}_N^\prime}(\infty) \leq (\geq) 0$;
    \item[(iv)] $H_{\mathcal{A}_N, \mathcal{B}_N}(\infty) \leq (\geq) 0$ and $H_{\mathcal{A}_N^\prime, \mathcal{B}_N^\prime}(\infty) > (<) 0$.
    \end{itemize}
\item[(C3)] there exists $x_2,x_3 \in (0,\infty)$ such that the function $x\mapsto \mathcal{A}_N(x)/\mathcal{B}_N(x)$ is increasing (decreasing) on $(0,x_2] \cup [x_3,\infty)$ and decreasing (increasing) on $[x_2,x_3]$ if:
    \begin{itemize}
    \item[(v)] $H_{\mathcal{A}_N, \mathcal{B}_N}(\infty) \geq (\leq) 0$, $H_{\mathcal{A}_N^\prime, \mathcal{B}_N^\prime}(\infty) > (<) 0$, and there exists $x_0\in(0,\infty)$ such that $H_{\mathcal{A}_N, \mathcal{B}_N}(x_0) < (>) 0$.
    \end{itemize}
\end{itemize}
\end{mr}

As we known, functions that satisfy certain conditions can be expanded into a power series, the so-called Maclaurin series.
Taking $a_k=\frac{\mathcal{A}^{(k)}(0)}{k!}$ and $b_k=\frac{\mathcal{B}^{(k)}(0)}{k!}$, power series $\mathcal{A}(x)$ and $\mathcal{B}(x)$ reduces to the so-called Maclaurin series with the following corollary holds.

\begin{cor} \label{cor-Taylor}
Let functions $\mathcal{A}(x)$ and $\mathcal{B}(x)$ have arbitrary order derivatives at point $t=0$,
as well as $\mathcal{A}(x)=\sum_{k=0}^\infty \frac{\mathcal{A}^{(k)}(0)}{k!} t^{k}$ and $\mathcal{B}(x)=\sum_{k=0}^\infty \frac{\mathcal{B}^{(k)}(0)}{k!} t^{k}$ converge on $(-r,r)$ with $\mathcal{B}^{(k)}(0)>0$ for all $k\geq0$.
If there exists different integers $m_2>m_1\geq1$ such that the sequence $\left \{\frac{\mathcal{A}^{(k)}(0)}{\mathcal{B}^{(k)}(0)} \right\}$ is increasing (decreasing) for all $k \in [0,m_1]_\mathbb{N} \cup [m_2,\infty)_\mathbb{N}$ and decreasing (increasing) for all $k \in [m_1,m_2]_\mathbb{N}$, then the conclusions (C1), (C2), and (C3) in Theorem \ref{DW21} hold.
\end{cor}

Now we consider some special cases by taking different function $\mathcal{B}$.
Taking $\mathcal{B}(x)= e^x = \sum_{k=0}^\infty \frac{x^k}{k!} $ in Monotonicity rule \ref{DW21}, then we have the following corollary.

\begin{cor}
Let real power series $\mathcal{A}(x)=\sum_{k=0}^\infty a_k x^k$ converge on $(0,\infty)$. If there exists different integers $m_2>m_1\geq1$ such that the sequence $\{k! a_k\}$ is increasing (decreasing) on $k \in [0,m_1]_\mathbb{N} \cup [m_2,\infty)_\mathbb{N}$ and decreasing (increasing) on $k \in [m_1,m_2]_\mathbb{N}$, then
\begin{itemize}
\item[(C1)] the function $x\mapsto e^{-x} \mathcal{A}(x)$ is increasing (decreasing) if one of the following conditions holds:
    \begin{itemize}
    \item[(i)] $\mathcal{A}^\prime(\infty) \geq (\leq) \mathcal{A}(\infty)$ and $\mathcal{A}^{\prime\prime}(\infty) \leq (\geq)  \mathcal{A}^\prime(\infty)$;
    \item[(ii)] $\mathcal{A}^\prime(\infty) > (<) \mathcal{A}(\infty)$, $\mathcal{A}^{\prime\prime}(\infty) > (<)  \mathcal{A}^\prime(\infty)$, and $\mathcal{A}^\prime(x) - \mathcal{A}(x) \geq 0$ for all $x\in(0,\infty)$.
    \end{itemize}
\item[(C2)] there exists $x_1 \in (0,\infty)$ such that the function $x\mapsto e^{-x} \mathcal{A}(x)$ is increasing (decreasing) on $(0,x_1]$ and decreasing (increasing) on $[x_1,\infty)$ if one of the following conditions holds:
    \begin{itemize}
    \item[(iii)] $\mathcal{A}^\prime(\infty) < (>) \mathcal{A}(\infty)$ and $\mathcal{A}^{\prime\prime}(\infty) \leq (\geq)  \mathcal{A}^\prime(\infty)$;
    \item[(iv)] $\mathcal{A}^\prime(\infty) \leq (\geq) \mathcal{A}(\infty)$ and $\mathcal{A}^{\prime\prime}(\infty) > (<)  \mathcal{A}^\prime(\infty)$.
    \end{itemize}
\item[(C3)] there exists $x_2,x_3 \in (0,\infty)$ such that the function $x\mapsto e^{-x} \mathcal{A}(x)$ is increasing (decreasing) on $(0,x_2] \cup [x_3,\infty)$ and decreasing (increasing) on $[x_2,x_3]$ if:
    \begin{itemize}
    \item[(v)] $\mathcal{A}^\prime(\infty) \geq (\leq) \mathcal{A}(\infty)$, $\mathcal{A}^{\prime\prime}(\infty) > (<)  \mathcal{A}^\prime(\infty)$, and there exists $x_0\in(0,\infty)$ such that $\mathcal{A}^\prime(x_0) < (>) \mathcal{A}(x_0)$.
    \end{itemize}
\end{itemize}
\end{cor}

Taking $\mathcal{B}(x)= \frac1{(1-x)^d} = \sum_{k=0}^\infty \frac{(d)_k}{k!} x^k$, then we have the following corollary, where $d>0$ and $(a)_k=\frac{\Gamma(a+k)}{\Gamma(a)}$ is the Pochhammer symbol.
\begin{cor}
Let real power series $\mathcal{A}(x)=\sum_{k=0}^\infty a_k x^k$ converge on $(0,1)$ and $d>0$. If there exists different integers $m_2>m_1\geq1$ such that the sequence $\{k! a_k/ (d)_k \}$ is increasing (decreasing) on $k \in [0,m_1]_\mathbb{N} \cup [m_2,\infty)_\mathbb{N}$ and decreasing (increasing) on $k \in [m_1,m_2]_\mathbb{N}$, then
\begin{itemize}
\item[(C1)] the function $x\mapsto (1-x)^{-d} \mathcal{A}(x)$ is increasing (decreasing) if one of the following conditions holds:
    \begin{itemize}
    \item[(i)] $\lim_{x \to 1} \big( \mathcal{A}^\prime(x) \frac{1-x}{d} - \mathcal{A}(x) \big) \geq (\leq) 0$ and
        $\lim_{x \to 1} \big( \mathcal{A}^{\prime\prime}(x) \frac{1-x}{d+1} - \mathcal{A}^\prime(x) \big) \leq (\geq) 0$;
    \item[(ii)] $\lim_{x \to 1} \big( \mathcal{A}^\prime(x) \frac{1-x}{d} - \mathcal{A}(x) \big) > (<) 0$, $\lim_{x \to 1} \big( \mathcal{A}^{\prime\prime}(x) \frac{1-x}{d+1} - \mathcal{A}^\prime(x) \big) > (<)  0$, and $\mathcal{A}^\prime(x) \frac{1-x}{d} - \mathcal{A}(x) \geq 0$ for all $x\in(0,1)$.
    \end{itemize}
\item[(C2)] there exists $x_1 \in (0,1)$ such that the function $x\mapsto (1-x)^{-d} \mathcal{A}(x)$ is increasing (decreasing) on $(0,x_1]$ and decreasing (increasing) on $[x_1,1)$ if one of the following conditions holds:
    \begin{itemize}
    \item[(iii)] $\lim_{x \to 1} \big( \mathcal{A}^\prime(x) \frac{1-x}{d} - \mathcal{A}(x) \big) < (>) 0$ and $\lim_{x \to 1} \big( \mathcal{A}^{\prime\prime}(x) \frac{1-x}{d+1} - \mathcal{A}^\prime(x) \big) \leq (\geq) 0$;
    \item[(iv)] $\lim_{x \to 1} \big( \mathcal{A}^\prime(x) \frac{1-x}{d} - \mathcal{A}(x) \big) \leq (\geq) 0$ and $\lim_{x \to 1} \big( \mathcal{A}^{\prime\prime}(x) \frac{1-x}{d+1} - \mathcal{A}^\prime(x) \big) > (<) 0$.
    \end{itemize}
\item[(C3)] there exists $x_2,x_3 \in (0,1)$ such that the function $x\mapsto (1-x)^{-d} \mathcal{A}(x)$ is increasing (decreasing) on $(0,x_2] \cup [x_3,1)$ and decreasing (increasing) on $[x_2,x_3]$ if:
    \begin{itemize}
    \item[(v)] $\lim_{x \to 1} \big( \mathcal{A}^\prime(x) \frac{1-x}{d} - \mathcal{A}(x) \big) \geq (\leq) 0$, $\lim_{x \to 1} \big( \mathcal{A}^{\prime\prime}(x) \frac{1-x}{d+1} - \mathcal{A}^\prime(x) \big) > (<) 0$, and there exists $x_0\in(0,1)$ such that $\mathcal{A}^\prime(x_0) \frac{1-x_0}{d} - \mathcal{A}(x_0) < (>) 0$.
    \end{itemize}
\end{itemize}
\end{cor}

In particular, taking $d=1$, namely, $\mathcal{B}(x)= \frac1{1-x} = \sum_{k=0}^\infty x^k$, then we have the following corollary.

\begin{cor}
Let real power series $\mathcal{A}(x)=\sum_{k=0}^\infty a_k x^k$ converge on $(0,1)$. If there exists different integers $m_2>m_1\geq1$ such that the sequence $\{a_k\}$ is increasing (decreasing) on $k \in [0,m_1]_\mathbb{N} \cup [m_2,\infty)_\mathbb{N}$ and decreasing (increasing) on $k \in [m_1,m_2]_\mathbb{N}$, then
\begin{itemize}
\item[(C1)] the function $x\mapsto (1-x) \mathcal{A}(x)$ is increasing (decreasing) if one of the following conditions holds:
    \begin{itemize}
    \item[(i)] $\lim_{x \to 1} \big( \mathcal{A}^\prime(x) (1-x) - \mathcal{A}(x) \big) \geq (\leq) 0$ and
        $\lim_{x \to 1} \big( \mathcal{A}^{\prime\prime}(x) \frac{1-x}{2} - \mathcal{A}^\prime(x) \big) \leq (\geq) 0$;
    \item[(ii)] $\lim_{x \to 1} \big( \mathcal{A}^\prime(x) (1-x) - \mathcal{A}(x) \big) > (<) 0$, $\lim_{x \to 1} \big( \mathcal{A}^{\prime\prime}(x) \frac{1-x}{2} - \mathcal{A}^\prime(x) \big) > (<)  0$, and $\mathcal{A}^\prime(x) (1-x) - \mathcal{A}(x) \geq 0$ for all $x\in(0,1)$.
    \end{itemize}
\item[(C2)] there exists $x_1 \in (0,1)$ such that the function $x\mapsto (1-x) \mathcal{A}(x)$ is increasing (decreasing) on $(0,x_1]$ and decreasing (increasing) on $[x_1,1)$ if one of the following conditions holds:
    \begin{itemize}
    \item[(iii)] $\lim_{x \to 1} \big( \mathcal{A}^\prime(x) (1-x) - \mathcal{A}(x) \big) < (>) 0$ and $\lim_{x \to 1} \big( \mathcal{A}^{\prime\prime}(x) \frac{1-x}{2} - \mathcal{A}^\prime(x) \big) \leq (\geq) 0$;
    \item[(iv)] $\lim_{x \to 1} \big( \mathcal{A}^\prime(x) (1-x) - \mathcal{A}(x) \big) \leq (\geq) 0$ and $\lim_{x \to 1} \big( \mathcal{A}^{\prime\prime}(x) \frac{1-x}{2} - \mathcal{A}^\prime(x) \big) > (<) 0$.
    \end{itemize}
\item[(C3)] there exists $x_2,x_3 \in (0,1)$ such that the function $x\mapsto (1-x) \mathcal{A}(x)$ is increasing (decreasing) on $(0,x_2] \cup [x_3,1)$ and decreasing (increasing) on $[x_2,x_3]$ if:
    \begin{itemize}
    \item[(v)] $\lim_{x \to 1} \big( \mathcal{A}^\prime(x) (1-x) - \mathcal{A}(x) \big) \geq (\leq) 0$, $\lim_{x \to 1} \big( \mathcal{A}^{\prime\prime}(x) \frac{1-x}{2} - \mathcal{A}^\prime(x) \big) > (<) 0$, and there exists $x_0\in(0,1)$ such that $\mathcal{A}^\prime(x_0) (1-x_0) - \mathcal{A}(x_0) < (>) 0$.
    \end{itemize}
\end{itemize}
\end{cor}

Taking $\mathcal{B}(x)= -\ln(1-d x) = \sum_{k=1}^\infty \frac{d^{k}}{(k)!} x^{k}$, then we have the following corollary.

\begin{cor}
Let real power series $\mathcal{A}(x)=\sum_{k=0}^\infty a_k x^k$ converge on $(0,1/d)$, where $d>0$. If there exists different integers $m_2>m_1\geq1$ such that the sequence $\{k! a_k / d^k\}$ is increasing (decreasing) on $k \in [0,m_1]_\mathbb{N} \cup [m_2,\infty)_\mathbb{N}$ and decreasing (increasing) on $k \in [m_1,m_2]_\mathbb{N}$, then
\begin{itemize}
\item[(C1)] the function $x\mapsto -\frac{\mathcal{A}(x)-a_0}{\ln(1-dx)} $ is increasing (decreasing) if one of the following conditions holds:
    \begin{itemize}
    \item[(i)] $\lim_{x \to \frac1{d}} \big( -\mathcal{A}^\prime(x) \frac{(1-d x) \ln(1-d x)}{d} - \mathcal{A}(x) \big) \geq (\leq) 0$ and
        $\lim_{x \to \frac1{d}} \big( \mathcal{A}^{\prime\prime}(x) \frac{1-d x}{d} - \mathcal{A}^\prime(x) \big) \leq (\geq) 0$;
    \item[(ii)] $\lim_{x \to \frac1{d}} \big( -\mathcal{A}^\prime(x) \frac{(1-d x) \ln(1-d x)}{d} - \mathcal{A}(x) \big) > (<) 0$, $\lim_{x \to \frac1{d}} \big( \mathcal{A}^{\prime\prime}(x) \frac{1-d x}{d} - \mathcal{A}^\prime(x) \big) > (<)  0$, and $-\mathcal{A}^\prime(x) \frac{(1-d x) \ln(1-d x)}{d} - \mathcal{A}(x) \geq 0$ for all $x\in(0,\frac1{d})$.
    \end{itemize}
\item[(C2)] there exists $x_1 \in (0,\frac1{d})$ such that the function $x\mapsto -\frac{\mathcal{A}(x)-a_0}{\ln(1-dx)} $ is increasing (decreasing) on $(0,x_1]$ and decreasing (increasing) on $[x_1,\frac1{d})$ if one of the following conditions holds:
    \begin{itemize}
    \item[(iii)] $\lim_{x \to \frac1{d}} \big( -\mathcal{A}^\prime(x) \frac{(1-d x) \ln(1-d x)}{d} - \mathcal{A}(x) \big) < (>) 0$ and $\lim_{x \to \frac1{d}} \big( \mathcal{A}^{\prime\prime}(x) \frac{1-dx}{d} - \mathcal{A}^\prime(x) \big) \leq (\geq) 0$;
    \item[(iv)] $\lim_{x \to \frac1{d}} \big( -\mathcal{A}^\prime(x) \frac{(1-d x) \ln(1-d x)}{d} - \mathcal{A}(x) \big) \leq (\geq) 0$ and $\lim_{x \to \frac1{d}} \big( \mathcal{A}^{\prime\prime}(x) \frac{1-dx}{d} - \mathcal{A}^\prime(x) \big) > (<) 0$.
    \end{itemize}
\item[(C3)] there exists $x_2,x_3 \in (0,\frac1{d})$ such that the function $x\mapsto -\frac{\mathcal{A}(x)-a_0}{\ln(1-dx)} $ is increasing (decreasing) on $(0,x_2] \cup [x_3,\frac1{d})$ and decreasing (increasing) on $[x_2,x_3]$ if:
    \begin{itemize}
    \item[(v)] $\lim_{x \to \frac1{d}} \big( -\mathcal{A}^\prime(x) \frac{(1-d x) \ln(1-d x)}{d} - \mathcal{A}(x) \big) \geq (\leq) 0$, $\lim_{x \to \frac1{d}} \big( \mathcal{A}^{\prime\prime}(x) \frac{1-dx}{d} - \mathcal{A}^\prime(x) \big) > (<) 0$, and there exists $x_0\in(0,\frac1{d})$ such that $-\mathcal{A}^\prime(x_0) \frac{(1-d x_0) \ln(1-d x_0)}{d}  - \mathcal{A}(x_0) < (>) 0$.
    \end{itemize}
\end{itemize}
\end{cor}

Respectively taking $\mathcal{B}(x)= \sinh(d x) = \sum_{k=0}^\infty \frac{d^{2k+1}}{(2k+1)!} x^{2k+1}$ and $\mathcal{B}(x)= \cosh(d x) = \sum_{k=0}^\infty \frac{d^{2k}}{(2k)!} x^{2k}$, then we have the following two corollaries.

\begin{cor}
Let real power series $\mathcal{A}(x)=\sum_{k=0}^\infty a_k x^{2k+1}$ converge on $(0,\infty)$ and $d>0$. If there exists different integers $m_2>m_1\geq1$ such that the sequence $\{(2k+1)! a_k / d^{2k+1}\}$ is increasing (decreasing) on $k \in [0,m_1]_\mathbb{N} \cup [m_2,\infty)_\mathbb{N}$ and decreasing (increasing) on $k \in [m_1,m_2]_\mathbb{N}$, then
\begin{itemize}
\item[(C1)] the function $x\mapsto \frac{\mathcal{A}(x)}{\sinh (dx)}$ is increasing (decreasing) if one of the following conditions holds:
    \begin{itemize}
    \item[(i)] $\lim_{x \to \infty} \big( \mathcal{A}^\prime(x) \frac{\tanh (d x)}{d} - \mathcal{A}(x) \big) \geq (\leq) 0$ and
        $\lim_{x \to \infty} \big( \mathcal{A}^{\prime\prime}(x) \frac{\coth (d x)}{d} - \mathcal{A}^\prime(x) \big) \leq (\geq) 0$;
    \item[(ii)] $\lim_{x \to \infty} \big( \mathcal{A}^\prime(x) \frac{\tanh (d x)}{d} - \mathcal{A}(x) \big) > (<) 0$, $\lim_{x \to \infty} \big( \mathcal{A}^{\prime\prime}(x) \frac{\coth (d x)}{d} - \mathcal{A}^\prime(x) \big) > (<)  0$, and $\mathcal{A}^\prime(x) \frac{\tanh (d x)}{d} - \mathcal{A}(x) \geq 0$ for all $x\in(0,\infty)$.
    \end{itemize}
\item[(C2)] there exists $x_1 \in (0,\infty)$ such that the function $x\mapsto \frac{\mathcal{A}(x)}{\sinh (dx)}$ is increasing (decreasing) on $(0,x_1]$ and decreasing (increasing) on $[x_1,\infty)$ if one of the following conditions holds:
    \begin{itemize}
    \item[(iii)] $\lim_{x \to \infty} \big( \mathcal{A}^\prime(x) \frac{\tanh (d x)}{d} - \mathcal{A}(x) \big) < (>) 0$ and $\lim_{x \to \infty} \big( \mathcal{A}^{\prime\prime}(x) \frac{\coth (d x)}{d} - \mathcal{A}^\prime(x) \big) \leq (\geq) 0$;
    \item[(iv)] $\lim_{x \to \infty} \big( \mathcal{A}^\prime(x) \frac{\tanh (d x)}{d} - \mathcal{A}(x) \big) \leq (\geq) 0$ and $\lim_{x \to \infty} \big( \mathcal{A}^{\prime\prime}(x) \frac{\coth (d x)}{d} - \mathcal{A}^\prime(x) \big) > (<) 0$.
    \end{itemize}
\item[(C3)] there exists $x_2,x_3 \in (0,\infty)$ such that the function $x\mapsto \frac{\mathcal{A}(x)}{\sinh (dx)}$ is increasing (decreasing) on $(0,x_2] \cup [x_3,\infty)$ and decreasing (increasing) on $[x_2,x_3]$ if:
    \begin{itemize}
    \item[(v)] $\lim_{x \to \infty} \big( \mathcal{A}^\prime(x) \frac{\tanh (d x)}{d} - \mathcal{A}(x) \big) \geq (\leq) 0$, $\lim_{x \to \infty} \big( \mathcal{A}^{\prime\prime}(x) \frac{\coth (d x)}{d} - \mathcal{A}^\prime(x) \big) > (<) 0$, and there exists $x_0\in(0,\infty)$ such that $\mathcal{A}^\prime(x_0) \frac{\tanh (d x_0)}{d} - \mathcal{A}(x_0) < (>) 0$.
    \end{itemize}
\end{itemize}
\end{cor}

\begin{cor}
Let real power series $\mathcal{A}(x)=\sum_{k=0}^\infty a_k x^{2k}$ converge on $(0,\infty)$ and $d>0$. If there exists different integers $m_2>m_1\geq1$ such that the sequence $\{(2k)! a_k / d^{2k} \}$ is increasing (decreasing) on $k \in [0,m_1]_\mathbb{N} \cup [m_2,\infty)_\mathbb{N}$ and decreasing (increasing) on $k \in [m_1,m_2]_\mathbb{N}$, then
\begin{itemize}
\item[(C1)] the function $x\mapsto \frac{\mathcal{A}(x)}{\cosh (dx)}$ is increasing (decreasing) if one of the following conditions holds:
    \begin{itemize}
    \item[(i)] $\lim_{x \to \infty} \big( \mathcal{A}^\prime(x) \frac{\coth (d x)}{d} - \mathcal{A}(x) \big) \geq (\leq) 0$ and
        $\lim_{x \to \infty} \big( \mathcal{A}^{\prime\prime}(x) \frac{\tanh (d x)}{d} - \mathcal{A}^\prime(x) \big) \leq (\geq) 0$;
    \item[(ii)] $\lim_{x \to \infty} \big( \mathcal{A}^\prime(x) \frac{\coth (d x)}{d} - \mathcal{A}(x) \big) > (<) 0$, $\lim_{x \to \infty} \big( \mathcal{A}^{\prime\prime}(x) \frac{\tanh (d x)}{d} - \mathcal{A}^\prime(x) \big) > (<)  0$, and $\mathcal{A}^\prime(x) \frac{\coth (d x)}{d} - \mathcal{A}(x) \geq 0$ for all $x\in(0,\infty)$.
    \end{itemize}
\item[(C2)] there exists $x_1 \in (0,\infty)$ such that the function $x\mapsto \frac{\mathcal{A}(x)}{\cosh (dx)}$ is increasing (decreasing) on $(0,x_1]$ and decreasing (increasing) on $[x_1,\infty)$ if one of the following conditions holds:
    \begin{itemize}
    \item[(iii)] $\lim_{x \to \infty} \big( \mathcal{A}^\prime(x) \frac{\coth (d x)}{d} - \mathcal{A}(x) \big) < (>) 0$ and $\lim_{x \to \infty} \big( \mathcal{A}^{\prime\prime}(x) \frac{\tanh (d x)}{d} - \mathcal{A}^\prime(x) \big) \leq (\geq) 0$;
    \item[(iv)] $\lim_{x \to \infty} \big( \mathcal{A}^\prime(x) \frac{\coth (d x)}{d} - \mathcal{A}(x) \big) \leq (\geq) 0$ and $\lim_{x \to \infty} \big( \mathcal{A}^{\prime\prime}(x) \frac{\tanh (d x)}{d} - \mathcal{A}^\prime(x) \big) > (<) 0$.
    \end{itemize}
\item[(C3)] there exists $x_2,x_3 \in (0,\infty)$ such that the function $x\mapsto \frac{\mathcal{A}(x)}{\cosh (dx)}$ is increasing (decreasing) on $(0,x_2] \cup [x_3,\infty)$ and decreasing (increasing) on $[x_2,x_3]$ if:
    \begin{itemize}
    \item[(v)] $\lim_{x \to \infty} \big( \mathcal{A}^\prime(x) \frac{\coth (d x)}{d} - \mathcal{A}(x) \big) \geq (\leq) 0$, $\lim_{x \to \infty} \big( \mathcal{A}^{\prime\prime}(x) \frac{\tanh (d x)}{d} - \mathcal{A}^\prime(x) \big) > (<) 0$, and there exists $x_0\in(0,\infty)$ such that $\mathcal{A}^\prime(x_0) \frac{\coth (d x_0)}{d} - \mathcal{A}(x_0) < (>) 0$.
    \end{itemize}
\end{itemize}
\end{cor}

\subsection{Other type monotonicity rules}
The following monotonicity rule shows that the monotonicity of $x\mapsto \mathcal{A}(x) / \mathcal{B}(x)$ near point $x=0$ is determined by the monotonicity of $\{a_k/b_k\}$ near point $k=0$, which named local monotonicity rule in this paper.

\begin{mr} \label{mr-add-1}
Let real power series $\mathcal{A}(x)=\sum_{k=0}^\infty a_k x^k$ and $\mathcal{B}(x)=\sum_{k=0}^\infty b_k x^k$ converge on $(0,r)$ with $b_k>0$. If there exists integer $m \geq 1$ such that the sequence $\{a_k/b_k\}$ is strictly increasing (decreasing) for all $0\leq k\leq m$, then there exists $x_0 \in (0,r]$ such that the function $x\mapsto \mathcal{A}(x) / \mathcal{B}(x)$ is strictly increasing (decreasing) on $(0,x_0)$.
\end{mr}

\begin{proof}
Noting that
\begin{equation*}
H_{\mathcal{A}, \mathcal{B}}(0^+) = \frac{\mathcal{A}^\prime(0^+)}{\mathcal{B}^\prime(0^+)}\mathcal{B}(0^+)- \mathcal{A}(0^+) = b_0 \Big( \frac{a_1}{b_1} - \frac{a_0}{b_0} \Big) > 0,
\end{equation*}
by using identity $(\mathcal{A}/\mathcal{B})^\prime = \mathcal{B}^\prime/ \mathcal{B}^2 H_{\mathcal{A},\mathcal{B}}$, we obtain that there exists $x_0 \in (0,r]$ such that $H_{\mathcal{A},\mathcal{B}}(x) > 0$ for all $x\in (0,x_0)$ and the function $x\mapsto \mathcal{A}(x) / \mathcal{B}(x)$ is strictly increasing (decreasing) on $(0,x_0)$.
\end{proof}

Let the monotonicity of $\{a_k/b_k\}$ changes $n(n\geq0)$ times and the the monotonicity of the response function$x\mapsto \mathcal{A}(x) / \mathcal{B}(x)$ changes $\tau_{\mathcal{A},\mathcal{B}}(n)$ times. The following monotonicity rule shows $\tau_{\mathcal{A},\mathcal{B}}(n) \leq n$ for all power series $\mathcal{A},\mathcal{B}$.

\begin{mr}
Let real power series $\mathcal{A}(x)=\sum_{k=0}^\infty a_k x^k$ and $\mathcal{B}(x)=\sum_{k=0}^\infty b_k x^k$ converge on $(0,r)$ with $b_k>0$. If the monotonicity of $\{a_k/b_k\}$ changes $n(n\geq0)$ times, then the monotonicity of $x\mapsto \mathcal{A}(x) / \mathcal{B}(x)$ changes no more than $n(n\geq0)$ times.
\end{mr}

\begin{proof}
We complete the proof by using mathematical induction. According to Monotonicity rule \ref{mr-SER-1}, Monotonicity rule \ref{mr-SER-5}, and Monotonicity rule \ref{DW21}, the conclusion is true for $n=0,1,2$, respectively. Suppose it is true for $n=q$, now we consider $n=q+1$.

In this case, the monotonicity of the sequence $\{a_k/b_k\}$ changes $q+1$ times. Without loss of generality, we suppose that there exist integer $m\geq1$, which can't be bigger, such that the sequence $\{a_k/b_k\}$ is increasing for all $0\leq k\leq m$.
Thus, the monotonicity of the sequence $\{a_{k+1}/b_{k+1}\}$ changes $q$ times if $m=1$ and changes $q+1$ times if $m\geq 2$.

(I) When $m=1$, the monotonicity of $x\mapsto \mathcal{A}^\prime(x) / \mathcal{B}^\prime(x)$ changes $q$ times with decreasing first, same as $H_{\mathcal{A},\mathcal{B}}$ due to identity $H_{\mathcal{A},\mathcal{B}}^\prime = (\mathcal{A}^\prime / \mathcal{B}^\prime )^\prime \mathcal{B}$. According to $H_{\mathcal{A},\mathcal{B}}(0^+)\geq0$, we obtain that the monotonicity of $x\mapsto \mathcal{A}(x) / \mathcal{B}(x)$ changes no more than $q+1$ times.

(II) Suppose it holds for $m=p$, we consider $m=p+1$. The monotonicity of $x\mapsto \mathcal{A}^\prime(x) / \mathcal{B}^\prime(x)$ changes $q+1$ times with increasing first, same as $H_{\mathcal{A},\mathcal{B}}$. According to $H_{\mathcal{A},\mathcal{B}}(0^+)\geq0$, we obtain that the monotonicity of $x\mapsto \mathcal{A}(x) / \mathcal{B}(x)$ changes no more than $q+1$ times.

Thus, we complete the proof.
\end{proof}

\begin{cor}
Let real power series $\mathcal{A}(x)=\sum_{k=0}^\infty a_k x^k$ and $\mathcal{B}(x)=\sum_{k=0}^\infty b_k x^k$ converge on $(0,r)$ with $b_k>0$ as well as $l\in\mathbb{N}$. If the monotonicity of $\{a_k/b_k\}$ changes $n(n\geq0)$ times, then the monotonicity of $x\mapsto \mathcal{A}^{(l)}(x) / \mathcal{B}^{(l)}(x)$ changes no more than $n$ times.
\end{cor}

\begin{proof}
Clearly, the monotonicity of $\{a_{k+l}/b_{k+l}\}$ changes no more than $n$ times. Thus, the monotonicity of $x\mapsto \mathcal{A}^{(l)}(x) / \mathcal{B}^{(l)}(x)$ changes no more than $n$ times.
\end{proof}

\end{document}